\documentclass[reqno,9 pt]{amsart}
\usepackage{amsmath,amsxtra,amssymb,latexsym, amscd,amsthm}
\usepackage[unicode]{hyperref}
\usepackage{array,tabularx,longtable,multicol,indentfirst,fancybox,color}
\usepackage{graphicx}
\usepackage{multicol}
\usepackage{mathrsfs}
\usepackage{enumerate}
\usepackage{tikz}
\usepackage{tikz-cd} 
\usepackage{tikz,tkz-tab}
\usetikzlibrary{calc, angles, quotes, intersections, positioning, patterns}
 
\advance\hoffset-1truecm\relax
\newtheorem{proposition}{Proposition}[section]
\newtheorem{theorem}{Theorem}[section]
\newtheorem{definition}{Definition}
\newtheorem{lemma}{Lemma}[section]
\newtheorem{corollary}{Corollary}[section]
\newtheorem{remark}{Remark}

\numberwithin{equation}{section}

\def\R{\mathbb{R}}

\def\lv{\lVert}

\numberwithin{equation}{section}

\begin{document}
	\title[Trigonometric weighted generalized convolution operator associated with  Fourier cosine-sine and  Kontorovich-Lebedev]{Trigonometric weighted generalized convolution operator associated with  Fourier cosine-sine and  Kontorovich-Lebedev transformations}
	
	\author[Trinh Tuan \& Nguyen Thanh Hong]{Trinh Tuan$^{1*}$  \& Nguyen Thanh Hong$^2$}
	\thanks{$^*$Corresponding author}
	\thanks{Received: 24 August 2023} 
	\thanks{Accepted: 4 April 2024 by \textit{Integral Transforms and Special Functions}}
	\thanks{Published online: 11 April 2024. \url{https://doi.org/10.1080/10652469.2024.2341400
	}}
	
	\maketitle	
\begin{center}
	$^1$Department of Mathematics, Faculty of Sciences, Electric Power University,\\ 235-Hoang Quoc Viet Rd., Bac Tu Liem Dist., Hanoi, Vietnam.\\
	E-mail: \texttt{tuantrinhpsac@yahoo.com}
	\\$^2$School for Gifted Students, Hanoi National University of Education,\\ 136-Xuan Thuy Rd., Cau Giay Dist., Hanoi, Vietnam. \\
	E-mail: \texttt{hongdhsp1@yahoo.com}
\end{center}	
\begin{abstract}
The main objective of this work is to introduce the generalized convolution with trigonometric weighted $\gamma=\sin y$ involving the Fourier cosine-sine and Kontorovich-Lebedev transforms, and to study its fundamental results. We establish the boundedness properties in a two-parametric family of Lebesgue spaces for this convolution operator. Norm estimation in the weighted $ L_p$ space is obtained and applications of the corresponding class of convolution integro-differential equations are discussed. The conditions for the solvability of these equations in $L_1$ space are also founded.
	\vskip 0.3cm
	\noindent\textbf{Keywords} Kontorovich–Lebedev transforms. Fourier transforms.  Convolution. Trigonometric weighted.
	\vskip 0.3cm
	
	\noindent \textbf{AMS Classifications}  42A38, 44A15,  44A35, 45A05. 
\end{abstract}	

\section{Introduction}
The theory of convolution of integral transform always remained an intensive area of exploration for researchers working in the area of mathematics, engineering, and physics. The essence of an integral transform lies in the property of the kernel function involved in it. When there are more than two transformations acting on the convolution we get the concept of generalized convolution. Studies of the various properties involved in an integral transformation are essentially kernel-dependent with respect to the integral. Therefore, in the process of development, various integral transformations have been generalized and further investigations have been carried out like the Fourier transform, Kontorovich–Lebedev transform, etc... In more detail, the reader can be referred to in \cite{Yaku94Luchko,Yakubovich2009Britvina,yakubovich2010Britvina}. These results lead to many approaches to solving differential equations, integral equations, image processing, parabolic type equations, and the boundedness of one-dimensional acoustic fields \cite{hoangtuan2017thaovkt,tuan2016generalized}.

Being directly influenced by the above and derived from  Yakubovich-Britvina's results in \cite{Yakubovich2009Britvina,yakubovich2010Britvina} was also the driving motivation for this work. We propose the definition of a generalized convolution $\big(f \overset{\gamma}{\underset{F_c,F_s,K}{\ast}} g\big)$ with trigonometric weighted $\gamma=\sin y$ and discuss some of its basic properties, which are also the main contributions of this article. We briefly recall some notions and results coming from  \cite{bateman1954,Sneddon1972,Titchmarsh1986,Yakubovich1996index,Prudnikov1986Marichev}. The Fourier cosine and Fourier sine transforms of the function $f$, denoted by  $(F_c)$ and  $(F_s)$ respectively, are defined by the integral formulas as follows
\begin{equation}\label{eq1.2}
(F_cf)(y):= \sqrt{\frac{2}{\pi}}\int_{\R_+} \cos xy f(x) dx,\ y>0, \ \text{and} \ (F_sf)(y):= \sqrt{\frac{2}{\pi}}\int_{\R_+} \sin xy f(x) dx,\ y>0.
\end{equation} 
The Kontorovich–Lebedev transform (often abbreviated as KL-transform) was introduced for the first time in \cite{KLtransform1938} to solve certain
boundary-value problems of mathematical physics. It arises naturally when the method of separation of variables is used to solve boundary-value problems formulated in terms of cylindrical
coordinate systems. Within the framework of this article, follow \cite{Yakubovich1996index} KL-transform denoted by  $(K)$, is defined by
$K[f](y):= \int_{\R_+} K_{iy}(x) f(x) dx,$
where $K_{iy} (x)$ is the modified Bessel's function \cite{Prudnikov1986Marichev} can be represented by
$
K_{iy}(x)=\int_{\R_+} e^{-x\cosh u} \cos yu du,$ $x>0$. 	The convolution related to  Fourier transform was introduced \cite{bateman1954}
 defined by 
 \begin{equation}\label{eq1.5}
 (f \underset{F}{\ast} g)(x) := \frac{1}{\sqrt{2\pi}}\int_{\R} f(x-u) g(u) du,\ x\in\mathbb{R}.
 \end{equation}
  Let $f,g$ be the functions belonging to $L_1 (\R)$, then the Fourier convolution is $ (f\underset{F}{\ast}g)\in L_1 (\R)$ (see \cite{Titchmarsh1986}). Furthermore, the factorization equality $
  F(f\underset{F}{\ast} g)(y) = (Ff)(y) (Fg)(y)
  $ is valid for any $y \in \R$.
According to \cite{Sneddon1972}, we consider the convolution of two functions $f$ and $g$ for the Fourier cosine transform defined by 
\begin{equation}\label{eq1.6}
(f \underset{F_c}{\ast} g)(x):= \frac{1}{\sqrt{2\pi}}\int_{\R_+} f(u)[g(|x-u|) + g(x+u)] du,\ x>0.
\end{equation}
If $f,g$  belong to $L_1 (\R_+)$, then  $(f \underset{F_c}{*}g) \in L_1 (\R_+)$ (refer\cite{Titchmarsh1986}), and the following factorization equality holds
\begin{equation}\label{eq1.4}
F_c (f \underset{F_c}{*}g)(y)=(F_c f)(y)(F_c g)(y), \ y>0.
\end{equation}
In addition, for any $f,g \in L_1 (\R_+)$ we obtain the $L_1$-norm estimation of Fourier cosine convolution as follow 
	\begin{equation}\label{eq1.8}
	\| f\underset{F_c}{\ast} g\|_{L_1(\mathbb{R}_+)} \leq 2\sqrt{\frac{2}{\pi}}\|f\|_{L_1(\mathbb{R}_+)} \|g\|_{L_1(\mathbb{R}_+)}.
	\end{equation}

\noindent Throughout this paper,  we shall make frequent use of weighted Lebesgue spaces $L_p(\mathbb{R}_+, \rho(x)), 1\leq p\leq \infty$ with respect to a positive measure $\rho(x)dx$ equipped with the norm $\|f\|_{L_p(\mathbb{R}_+;\rho)}=\big(\int_{\R_+} |f(x)|^p \rho(x) dx\big)^{1\textfractionsolidus p}$ is finite. In case $\rho(x)= K_0(\beta x)x^\alpha$ with $\alpha\in \mathbb{R}$, $0<\beta\leq 1$, we have a two-parametric family of Lebesgue spaces $L_p^{\alpha,\beta}(\mathbb{R}_+) \equiv L_p (\mathbb{R}_+;K_0(\beta x)x^\alpha)$ defined by $
L_p^{\alpha,\beta}(\mathbb{R}_+) := \left\{f(x): \int_{\R_+} |f(x)|^p K_0(\beta x)x^\alpha dx<\infty\right\},
$ and normed by  $
\|f\|_{L_p^{\alpha,\beta}(\mathbb{R}_+)} = \big(\int_{\R_+} |f(x)|^p K_0(\beta x) x^\alpha dx\big)^{1\textfractionsolidus p}
$ is finite \cite{yakubovich2010Britvina}.  

This paper is divided into four sections and is organized as follows. Section 2 is devoted to presenting the concept of the trigonometric weighted generalized convolution. The structurally important properties of the operator associated with defined spaces are also clearly established. We show the existence of this operator on space $L_1(\mathbb{R}_+)$, simultaneously proving the factorization equality and Parseval equality. Section 3 consists of two subsections containing the most important results of this article. In subsection 3.1, we give another version of the Young-type theorem for $\big(f \overset{\gamma}{\underset{F_c,F_s,K}{\ast}} g\big)$ and prove that it is a bounded function on $L_{\infty} (\R_+)$. The general formulation of the Young-type inequality for generalized convolution \eqref{eq2.1} is investigated through Riesz's representation theorem. The above results are key materials to find out the boundedness of convolution \eqref{eq2.1} on a two-parametric family of Lebesgue spaces $L_s^{\gamma_1,\gamma_2}(\mathbb{R}_+)$. A sharp upper bound in this estimate is specifically expressed through Euler's gamma-function. In subsection 3.2, using H\"older's inequality, Fubini's theorem, and modified Bessel's function, we establish norm inequalities in the weighted $L_p$ spaces. Section 4 displays several applications of the constructed generalized convolution for the solvability of classes of convolution integro-differential equations. Namely, by using obtained results together with the help of Wiener-L\'evy's theorem, we provide the conditions for the solvability of the first and second kinds of integro-differential equations involving operator \eqref{eq2.1} and obtain explicit $L_1$-solutions. 
\section{Structure of trigonometric weighted generalized convolution}
\begin{definition}\label{def2.1}
	The generalized convolution operator for the Fourier cosine-sine, and Kontorovich-Lebedev integral transforms with trigonometric weighted $\gamma(y)=\sin y$ of two functions $f$, $g$ is denoted by $(f \overset{\gamma}{\underset{F_c,F_s,K}{\ast}} g)$ and defined by
	\begin{equation}\label{eq2.1}
	(f \overset{\gamma}{\underset{F_c,F_s,K}{\ast}} g)(x):=\frac14\int_{\mathbb{R}^2_+ } \varphi(x,u,v) f(u) g(v) du dv,\quad x>0,
	\end{equation}
	with the kernel function
	\begin{equation}\label{eq2.2}
	\varphi(x,u,v)=e^{-v \cosh(x+u-1)} + e^{-v\cosh(x-u+1)} - e^{-v\cosh(x+u+1)} - e^{-v\cosh(x-u-1)}.
	\end{equation}
\end{definition}
\begin{theorem}\label{thm2.1}
	Suppose that $f$ is an arbitrary function in $ L_1(\mathbb{R}_+)$ and $g\in L_1^{0,\beta}(\mathbb{R}_+)$, we have\\
	\textbf{i)} For  $\beta \in (0,1]$, then the convolution \eqref{eq2.1} is well-defined  for all $x >0$ as a continuous function and belongs to $L_1(\mathbb{R}_+)$. Besides, we obtain the $L_1$-norm estimation  as follows 
	\begin{equation}\label{eq2.3}
	\|f \overset{\gamma}{\underset{F_c,F_s,K}{\ast}} g\|_{L_1(\mathbb{R}_+)} \leq 2\|f\|_{L_1(\mathbb{R}_+)} \|g\|_{L_1^{0,\beta}(\mathbb{R}_+)}.
	\end{equation}
	\textbf{ii)} In case $0<\beta<1$, for all $x>0$, then convolution \eqref{eq2.1} satisfies the generalized Parseval type equality 
	\begin{equation}\label{eq2.4}
	(f \overset{\gamma}{\underset{F_c,F_s,K}{\ast}} g)(x)=\sqrt{\frac{2}{\pi}}\int_{\R_+} (F_sf)(y) K[g](y) \sin y \cos(xy) dy,
	\end{equation}
	and the following factorization property is valid
	\begin{equation}\label{eq2.5}
	F_c(f\overset{\gamma}{\underset{F_c,F_s,K}{\ast}}g)(y)=\sin y(F_sf)(y) K[g](y).
	\end{equation}		
	Furthermore $(f \overset{\gamma}{\underset{F_c,F_s,K}{\ast}} g)\in C_0(\mathbb{R}_+)$, where $C_0(\mathbb{R}_+)$ is the space of bounded continuous functions vanishing at infinity.

\end{theorem}

\begin{proof}
	To prove assertion \textbf{\textit{i)}}, we need to show $\int_{\R_+} \big|(f \overset{\gamma}{\underset{F_c,F_s,K}{\ast}} g)(x)\big| dx$ is finite.
	Indeed, for any $v,x>0$ we have $\int_{\R_+} e^{-v\cosh (x+u-1)} du = \int_{x-1}^\infty e^{-v\cosh t} dt \leq \int_{\R} e^{-v\cosh t}dt= 2K_0(v)\leq 2K_0(\beta v), \forall \beta \in (0,1]$.
	By definition \eqref{eq2.2}, we get estimate of the kernel as follows $\int_{\R_+} |\varphi(x,u,v)|du = \int_{\R_+} |\varphi(x,u,v)| dx \leq 8K_0(\beta v).$ 
	Coupling the above with  Fubini's theorem, we obtain 
	\begin{align*}
	\int_{\R_+} \big|(f \overset{\gamma}{\underset{F_c,F_s,K}{\ast}} g)(x)\big| dx&\leq \frac14\int_{\mathbb{R}_+^3} |\varphi(x,u,v)| |f(u)| |g(v)| du dv dx\\
	&=\frac14\int_{\R_+} |f(u)|\left[\int_{\R_+} |g(v)|\bigg(\int_{\R_+}|\varphi(x,u,v)|dx\bigg)\right]dv du\\
	&\leq \frac14\int_{\R_+} |f(u)| \left(\int_{\R_+} 8K_0(\beta v)|g(v)| dv\right)du\\
	&=2\left(\int_{\R_+} |f(u)| du\right)\left(\int_{\R_+} K_0(\beta v)|g(v)| dv\right)
	=2\|f\|_{L_1(\mathbb{R}_+)} \|g\|_{L_1^{0,\beta}(\mathbb{R}_+)}.
	\end{align*}
	Thus $\|g\|_{L_1^{0,\beta}(\mathbb{R}_+)}= \big(\int_{\R_+} |g(x)| K_0(\beta x) dx\big)
	$ is finite, then $\int_{\R_+} \big|(f \overset{\gamma}{\underset{F_c,F_s,K}{\ast}} g)(x)\big| dx$ is finite for
	almost all $x > 0$, and it implies that
	$(f \overset{\gamma}{\underset{F_c,F_s,K}{\ast}} g)$ belongs to $L_1(\mathbb{R}_+)$ and we derive the estimate \eqref{eq2.3}.\\
	 \textbf{\textit{ii)}}
	Applying formula $2.16.48.19$ in \cite{Prudnikov1986Marichev}, we infer that
	$
	\int_{\R_+} \cos (ty) K_{iy}(v) dy =\frac{\pi}2e^{-v\cosh t}.
	$ Therefore
\begin{equation}\label{eq2.6}
\begin{aligned}
(f\overset{\gamma}{\underset{F_c,F_s,K}{\ast}} g) (x)&=\frac14\int_{\mathbb{R}_+^3} \frac2\pi K_{iy}(v)[\cos y(x+u-1) + \cos y(x-u+1)\\
&- \cos y(x+u+1) - \cos y(x-u-1)] f(u) g(v) du dv dy\\
&=\frac2\pi \int_{\mathbb{R}_+^3} K_{iy}(v) f(u) g(v) \sin y \cos(yx) \sin(yu)\, du dv dy.
\end{aligned}
\end{equation}
According to \cite{Yakubovich1996index}, we obtain  estimation $|K_{iy}(v)|\leq e^{-y\arccos\beta}K_0(\beta v)$ that holds true for all $y>0$, $v>0$, with $0<\beta< 1$. For any $f\in L_1(\mathbb{R}_+)$ and $g\in L_1^{0,\beta}(\mathbb{R}_+)$, we have
	$$\begin{aligned}
	&\int_{\mathbb{R}_+^3} |K_{iy}(v) f(u) g(v) \sin y \cos (yx) \sin (yu)|du dv dy
	\leq \int_{\mathbb{R}_+^3} |K_{iy}(v)| |f(u)| |g(v)| du dv dy\\
	&\leq \int_{\mathbb{R}_+^3} e^{-y\arccos \beta}K_0(\beta v) |f(u)| |g(v)| du dv dy
	=\left(\int_{\R_+} e^{-y\arccos \beta}dy\right)\left(\int_{\R_+} |f(u)| du\right)\left(\int_{\R_+} K_0(\beta v)|g(v)| dv\right)\\
	&=\frac{1}{\arccos\beta}\|f\|_{L_1(\mathbb{R}_+)}\|g\|_{L_1^{0,\beta}(\mathbb{R}_+)}<\infty.
	\end{aligned}$$
	This implies that integral \eqref{eq2.6} is absolutely convergent, by using Fubini's theorem we obtain
	\begin{align*}
	(f \overset{\gamma}{\underset{F_c,F_s,K}{\ast}} g)(x) &= \sqrt{\frac{2}{\pi}}\int_{\R_+}\left\{\left(\sqrt{\frac{2}{\pi}}\int_{\R_+} f(u)\sin (yu) du\right)
\left(\int_{\R_+} K_{iy}(v) g(v) dv\right)\sin y\cos (xy)\right\} dy.
	\end{align*}
Combining \eqref{eq1.2},\eqref{eq1.4}, we deduce that
	$
	(f \overset{\gamma}{\underset{F_c,F_s,K}{\ast}} g)(x) = \sqrt{\frac{2}{\pi}}\int_{\R_+}\sin y(F_sf)(y) K[g](y) \cos (xy) dy.
	$
	Applying the Fourier cosine transform $(F_c)$  on both sides of the Parseval equality, we derive the factorization equality \eqref{eq2.5}. The Riemann-Lebesgue's theorem in \cite{Sogge1993fourier} states that \textquotedblleft If $f\in L_1(\R^n)$, then $(Ff)(y) \rightarrow 0$ as  $|y|$ tends to $\infty$, and, hence $(Ff)(y)\in C_0 (\R^n)$\textquotedblright. This is still true for the  Fourier sine transform  $(F_s)$ on $\R_+$ (refer \cite{Titchmarsh1986}), implying that if $f$ belongs to $L_1(\mathbb{R}_+)$, then $(F_s f)\in C_0(\mathbb{R}_+)$ and $|(F_sf)(y)|\leq \sqrt{\frac{2}{\pi}}\|f\|_{L_1(\mathbb{R}_+)}$. Follow \cite{Yakubovich1996index}, we have an estimate $$|K[g](y)| \leq \int_{\R_+} |K_{iy}(x)| |g(x)| dx \leq \int_{\R_+} e^{-|y|\arccos\beta}K_0(\beta x) g(x) dx
	\leq e^{-|y|\arccos\beta}\ \|g\|_{L_1^{0,\beta}(\mathbb{R}_+)}<\infty.$$
	Therefore, $(F_sf)(y) K[g](y)$ is a bounded function on $\R_+$.  Letting $x$ tends to $\infty$ in Parseval equality \eqref{eq2.4}, we obtain the final conclusion of this theorem obvious by virtue of Riemann-Lebesgue's theorem.
\end{proof}

\begin{remark}\label{rem1'}
\textup{In the case $f\in L_2(\mathbb{R}_+)$, $g\in L_2^{0,\beta}(\mathbb{R}_+)$ then the operator \eqref{eq2.1} satisfies the equalities \eqref{eq2.4} and \eqref{eq2.5}, where the integrals are understood in mean-square
convergence sense. The integral converges uniformly on interval $[0,N]$, therefore, the limit is
	$$
(f \overset{\gamma}{\underset{F_c,F_s,K}{\ast}} g)(x)=\sqrt{\frac{2}{\pi}}\lim_{N\to \infty}\int_{0}^N (F_sf)(y) K[g](y) \sin y \cos(xy) dy,\quad x>0.
$$
 Here, we define the cosine and sine Fourier transforms in the mean-square convergence sense, namely
$$(F_{\left\{\substack{c\\s}\right\}}f)(y):= \sqrt{\frac{2}{\pi}}\lim\limits_{N\to \infty}\int_{0}^N \left\{\substack{\cos (xy)\\ \sin (xy)}\right\} f(x) dx,$$
and Plancherel’s theorem in \cite{Sneddon1972} said that $F_c ,F_s : L_2 (R_+) \longrightarrow L_2 (R_+)$ are
isometric isomorphisms mappings with Parseval’s equalities $\|F_{\left\{{ }_s^c\right\}} f\|_{L_2(\mathbb{R}_{+} )}=\|f\|_{L_2(\mathbb{R}_{+} )}.$ According to \cite{Yakubovich1996index}, we know that KL-transform is an isometric isomorphism mapping $K_{i y}$ : $L_2(\mathbb{R}_{+} ; x dx) \longrightarrow L_2(\mathbb{R}_{+} ; x \sinh \pi xdx)$, where integral $\int_{\R_+} K_{iy}(x) g(x) dx$ does not exist in Lebesgue's sense and therefore we understand it in the form 
$$
	K[g](y):= \lim_{N\to \infty}\int_{\frac1N}^\infty K_{iy}(x) g(x) dx.
$$}
\end{remark}

\section{Boundedness on a two-parametric family of Lebesgue spaces and estimation on weighted space}

\subsection{Boundedness on index spaces $L_s^{\gamma_1,\gamma_2}(\mathbb{R}_+)$}
The Young inequality for Fourier convolution operator was introduced in \cite{YoungWH1912}
$
\| f\underset{F}{\ast} g\|_{L_r(\mathbb{R}_+)} \leq \|f\|_{L_p(\mathbb{R})}\|g\|_{L_q(\mathbb{R})},
$ for any $f\in L_p(\mathbb{R})$, $g\in L_q(\mathbb{R})$,
where $p,q,r>1$ such that $\frac1p+\frac1q=1+\frac1r$. After that, this inequality  was extended by Adams--Fournier (Theo. 2.24 in \cite{AdamsFournier2003sobolev}) as follows 
$
\left|\int_{\mathbb{R}^n} (f \underset{F}{\ast} g) (x) \omega(x) dx\right| \leq \|f\|_{L_p(\mathbb{R}^n)} \|g\|_{L_q(\mathbb{R}^n)} \|\omega\|_{L_r(\mathbb{R}^n)}
$
here, $p, q, r>1$ such that $\frac1p+\frac1q+\frac1r=2$ and $f\in L_p(\mathbb{R}^n)$, $g\in L_q(\mathbb{R}^n)$, $\omega\in L_r(\mathbb{R}^n)$.
In this subsection, by using the technique as in \cite{AdamsFournier2003sobolev,Tuan2023VKT,tuan2020Ukrainian}, we introduce the Young type theorem for operator \eqref{eq2.1} and prove the boundedness in $L_r(\mathbb{R}_+)$ with $1<r<\infty$. The case $r=\infty$ is also discussed in detail. Moreover, with the help of the technique that is presented in the proof of Theorem  \ref{thm3.2}, we obtain a generalized result for the boundedness of operator \eqref{eq2.1} on the index space  $L_s^{\gamma_1,\gamma_2}(\mathbb{R}_+)$ with $s\geq 1$ and the parameters $\gamma_1>-1$, $\gamma_2 >0$.
\begin{theorem}[Young type theorem for convolution \eqref{eq2.1}]\label{thm3.1} Let $p,q$, and $r$ be real numbers in open interval $(1,\infty)$ such that $\frac{1}{p}+\frac{1}{q}+\frac{1}{r}=2$.	
	For any functions $f\in L_p(\mathbb{R_+})$, $g\in L^{0,\beta}_q(\mathbb{R_+})$, and $h\in L_r(\mathbb{R})$,  we obtain the following estimation with $0<\beta\leq 1$
	\begin{equation}\label{eq3.2}
	\left|\int_{\R_+}(f \overset{\gamma}{\underset{F_c,F_s,K}{\ast}}g) (x) h(x) dx\right| \leq 2^{\frac1q}\|f\|_{L_p(\mathbb{R}_+)} \|g\|_{L_q^{0,\beta}(\mathbb{R}_+)} \|h\|_{L_r(\mathbb{R}_+)}.
	\end{equation}
\end{theorem}
\begin{proof}
Let  $p_1,q_1,r_1$  be the conjugate exponentials of $p,q,r$, respectively. This means that  $\frac{1}{p}+\frac{1}{p_1}=\frac{1}{q}+\frac{1}{q_1}=\frac{1}{r}+\frac{1}{r_1}=1$, together with the assumption of theorem, we get the correlation between exponential numbers as follows
\begin{equation}\label{eq3.4}
\left\{\begin{array}{l}
\frac{1}{p_1}+\frac{1}{q_1}+\frac{1}{r_1}=1,\\
p\left(\frac{1}{q_1}+\frac{1}{r_1}\right)=q\left(\frac{1}{p_1}+\frac{1}{r_1}\right)=r\left(\frac{1}{p_1}+\frac{1}{q_1}\right)=1.
\end{array}\right.
\end{equation}
	For simplicity, we set	
	\begin{align*}
	T_1(x,u,v)&=|\varphi(x,u,v)|^{\frac1{p_1}} |g(v)|^{\frac{q}{p_1}} |h(x)|^{\frac{r}{p_1}}\in L_{p_1}(\mathbb{R}_+^3),\\
	T_2(x,u,v)&=|\varphi(x,u,v)|^{\frac{1}{q_1}}|f(u)|^{\frac{p}{q_1}}|h(x)|^{\frac{r}{q_1}}\in L_{q_1}(\mathbb{R}_+^3),\\
	T_3(x,u,v)&=|\varphi(x,u,v)|^{\frac{1}{r_1}}|f(u)|^{\frac{p}{r_1}}|g(v)|^{\frac{q}{r_1}}\in L_{r_1}(\mathbb{R}_+^3).
	\end{align*}
	Under the conditions \eqref{eq3.4}, we get $
	T_1(x,u,v) T_2(x,u,v) T_3(x,u,v)=|\varphi(x,u,v)| |f(u)| |g(v)| |h(x)|.
	$
	Therefore
	\begin{align*}
	I=\left|\int_{\R_+} (f \overset{\gamma}{\underset{F_c,F_s,K}{\ast}} g)(x) h(x) dx\right|\leq \frac14\int_{\mathbb{R}_+^3} T_1(x,u,v)  T_2(x,u,v) T_3(x,u,v) du dv dx.
	\end{align*}
	Moreover, since $\frac{1}{p_1}+\frac{1}{q_1}+\frac{1}{r_1}=1$, applying the H\"older inequality, we deduce that
	\begin{equation}\label{eq3.5}
	\begin{aligned}
	I&\leq \frac14\left\{\int_{\mathbb{R}_+^3} |T_1(x,u,v)|^{p_1} du dv dx\right\}^{\frac{1}{p_1}}\left\{\int_{\mathbb{R}_+^3}|T_2(x,u,v)|^{q_1} du dv dx\right\}^{\frac1{q_1}}
	 \left\{\int_{\mathbb{R}_+^3} |T_3(x,u,v)|^{r_1} du dv dx\right\}^{\frac{1}{r_1}}\\
	&=\frac14\|T_1\|_{L_{p_1}(\mathbb{R}_+^3)}  \|T_2\|_{L_{q_1}(\mathbb{R}_+^3)} \|T_3\|_{L_{r_1}(\mathbb{R}_+^3)}.
	\end{aligned}\end{equation}
Directly inferred from the formula \eqref{eq2.2} for any $\beta \in (0,1]$, we have 
\begin{equation}\label{danhgiaK0}
\int_{\R_+} |\varphi(x,u,v)|du\leq 8 K_0(v)\leq 8K_0(\beta v), \quad v>0	.
\end{equation}	
Based on the assumption of $f\in L_p(\mathbb{R}_+)$ with $g\in L_q^{0,\beta}(\mathbb{R}_+)$, $0<\beta\leq 1$ and $h\in L_r(\mathbb{R}_+)$, using Fubini's theorem and \eqref{danhgiaK0},  we obtain  $L_{p_1} (\R^3_+)$-norm estimation for the operator $T_1 $ as follows	
	\begin{align*}
	&\|T_1\|^{p_1}_{L_{p_1}(\mathbb{R}_+^3)}\\&=\int_{\mathbb{R}_+^3} \left\{|\varphi(x,u,v)|^{\frac{1}{p_1}} |g(v)|^{\frac{q}{p_1}} |h(x)|^{\frac{r}{p_1}}\right\}^{p_1} du dv dx
	= \int_{\R_+} |h(x)|^r\left\{\int_{\R_+} |g(v)|^q\left(\int_{\R_+} |\varphi(x,u,v)| du\right)dv\right\} dx\\
	&\leq \int_{\R_+} |h(x)|^r\left\{\int_{\R_+} 8|g(v)|^q K_0(v) dv\right\} dx
	\leq 8\int_{\R_+} |h(x)|^r\left(\int_{\R_+} K_0(\beta v) |g(v)|^q dv\right) dx\\
	&= 8\left(\int_{\R_+} |h(x)|^r dx\right) \left(\int_{\R_+} K_0(\beta v) |g(v)|^q dv\right)
	= 8\|h\|^r_{L_r(\mathbb{R}_+)} \|g\|^q_{L_q^{0,\beta}(\mathbb{R}_+)}.
	\end{align*}
	Therefore
	\begin{equation}\label{eq3.6}
	\|T_1\|_{L_{p_1}(\mathbb{R}_+^3)} \leq 8^{\frac1{p_1}} \|h\|^{\frac{r}{p_1}}_{L_r(\mathbb{R}_+)} \|g\|^{\frac{q}{p_1}}_{L_q^{0,\beta}(\mathbb{R}_+)}.
	\end{equation}
		Similar to what we did with the evaluation \eqref{eq3.6} of $T_1$, we also get the norm estimation of $T_3$ on $L_{r_1} (\R^3_+)$ as follows
	\begin{equation}\label{eq3.7}
	\|T_3\|_{L_{r_1}(\mathbb{R}_+^3)}\leq 8^{\frac{1}{r_1}}\|f\|^{\frac{p}{r_1}}_{L_p(\mathbb{R}_+)} \|g\|^{\frac{q}{r_1}}_{L_q^{0,\beta}(\mathbb{R}_+)}.
	\end{equation}
To give an estimate for operator $T_2$, it is easy to first see that $\int_{\R_+} |\varphi(x,u,v)|dv\leq 4\int_{\R_+} e^{-v} dv=4$. This means that
	$$\begin{aligned}
	\|T_2\|^{q_1}_{L_{q_1}(\mathbb{R}_+^3)}&=\int_{\mathbb{R}_+^3}|\varphi(x,u,v)||f(u)|^p|h(x)|^r du dv dx\\
	&=\left(\int_{\R_+} |\varphi(x,u,v)|dv\right)\left(\int_{\R_+} |f(u)|^p du\right)\left(\int_{\R_+} |h(x)|^rdx\right)
	\leq 4\|f\|^p_{L_p(\mathbb{R}_+)} \|h\|^r_{L_r(\mathbb{R}_+)}.
	\end{aligned}$$
	This yields
	\begin{equation}\label{eq3.8}
	\|T_2\|_{L_{q_1}(\mathbb{R}_+^3)} \leq 4^{\frac1{q_1}}\|f\|^{\frac{p}{q_1}}_{L_p(\mathbb{R}_+)} \|h\|^{\frac{r}{q_1}}_{L_r(\mathbb{R}_+)}.
	\end{equation}
	Coupling \eqref{eq3.6},\eqref{eq3.7}, and \eqref{eq3.8} we have
	\begin{equation}\label{eq3.9}
	\|T_1\|_{L_{p_1}(\mathbb{R}_+^3)} \|T_2\|_{L_{q_1}(\mathbb{R}_+^3)} \|T_3\|_{L_{r_1}(\mathbb{R}_+^3)} \leq 8\left(\frac12\right)^{\frac{1}{q_1}} \|f\|_{L_p(\mathbb{R}_+)} \|g\|_{L_q^{0,\beta}(\mathbb{R}_+)}\|h\|_{L_r(\mathbb{R}_+)}.
	\end{equation}
	Finally, combining \eqref{eq3.9} and \eqref{eq3.5}, we obtain estimation as in the conclusion of the theorem.
\end{proof}
The following Young-type inequality is a direct consequence of Theorem \ref{thm3.1}.
\begin{corollary}[Young type inequality for convolution \eqref{eq2.1}]\label{cor3.1}
	Let $p,q,r \in (1,\infty)$, satisfying  $\frac{1}{p}+\frac{1}{q}=1+\frac{1}{r}$. If $f\in L_p(\mathbb{R_+})$, $g\in L_q^{0,\beta} (\mathbb{R}_+)$ with $0<\beta\leq 1$, then the convolution \eqref{eq2.1} is well-defined and belongs to $L_r(\mathbb{R_+})$. Hence the following inequality holds
	
	\begin{equation}\label{eq3.10}
	\big\|f \overset{\gamma}{\underset{F_c,F_s,K}{\ast}} g\big\|_{L_r(\mathbb{R}_+)} \leq 2^{\frac1q} \|f\|_{L_p(\mathbb{R}_+)} \|g\|_{L_q^{0,\beta}(\mathbb{R}_+)}.
	\end{equation}
\end{corollary}
\begin{proof}
	Let $r_1$ be the conjugate exponent of $r$, i.e $\frac{1}{r}+\frac{1}{r_1}=1$. From the assumptions of Corollary \ref{cor3.1}, we have $\frac{1}{p}+\frac{1}{q}+\frac{1}{r_1}=2$,  which shows the numbers $p$, $q$, and $r_1$  satisfy the conditions of  Theorem \ref{thm2.1} (with role of $r$  being replaced by $r_1$). Therefore, if $f\in L_p(\mathbb{R}_+)$, $g\in L_q^{0,\beta}(\mathbb{R}_+)$ then the linear operator$$\mathscr{T}h:=\int_{\mathbb{R}_+} (f \overset{\gamma}{\underset{F_c,F_s,K}{\ast}} g)(x)\cdot h(x) dx$$ is bounded in $L_{r_1}(\R_+)$.
Consequently, by the Riesz's representation theorem \cite{Stein1971Weiss}, then generalized convolution $(f \overset{\gamma}{\underset{F_c,F_s,K}{\ast}} g)(x)$ belongs to $L_r(\mathbb{R_+})$.
To prove the inequality \eqref{eq3.10},  we choose the function $$h(x):=\mathrm{sign}\bigg\{ (f \overset{\gamma}{\underset{F_c,F_s,K}{\ast}} g)(x)\bigg\}^r \times \bigg\{ (f \overset{\gamma}{\underset{F_c,F_s,K}{\ast}} g)(x)\bigg\}^{\frac{r}{r_1}}.$$ 
Then $h \in L_{r_1} (\R_+)$, with the norm $
\|h\|_{L_{r_1}(\mathbb{R}_+)} = \big\|f \overset{\gamma}{\underset{F_c,F_s,K}{\ast}} g\big\|^{\frac{r}{r_1}}_{L_r(\mathbb{R}_+)}.
$ Applying inequality \eqref{eq3.2} to such  function $h(x)$, we get
\begin{align*}
\|f \overset{\gamma}{\underset{F_c,F_s,K}{\ast}} g\|^r_{L_r(\mathbb{R}_+)}&=\int_{\R_+} \big|(f \overset{\gamma}{\underset{F_c,F_s,K}{\ast}} g)(x)\big|^r dx
= \bigg|\int_{\R_+}(f \overset{\gamma}{\underset{F_c,F_s,K}{\ast}} g)(x) \cdot h(x) dx\bigg| \\
&\leq 2^{\frac1q}\|f\|_{L_p(\mathbb{R}_+)} \|g\|_{L_q^{0,\beta}(\mathbb{R}_+)} \|h\|_{L_{r_1}(\mathbb{R}_+)}
=2^{\frac1q}\|f\|_{L_p(\mathbb{R}_+)}\|g\|_{L_q^{0,\beta}(\mathbb{R}_+)}\|f \overset{\gamma}{\underset{F_c,F_s,K}{\ast}}g\|^{\frac{r}{r_1}}_{L_r(\mathbb{R}_+)},
\end{align*}or equivalent for any $0<\beta \leq 1$, then  $\|f \overset{\gamma}{\underset{F_c,F_s,K}{\ast}} g\|^{r-\frac{r}{r_1}}_{L_r(\mathbb{R}_+)} \leq 2^{\frac1q}\|f\|_{L_p(\mathbb{R}_+)}\|g\|_{L_q^{0,\beta}(\mathbb{R}_+)}$.
Since $r - \frac{r}{r_1} =1$, we arrive at the conclusion of the corollary.

\end{proof}

What about the case $r=\infty$? We consider the boundedness of operator \eqref{eq2.1} in the case $r=\infty$ via the following theorem.
\begin{theorem}\label{thm3.2}
	Suppose that $p,q>1$ and satisfy $\frac{1}{p}+\frac{1}{q}=1$. For any functions $f\in L_p(\mathbb{R}_+)$, $g\in L_q^{0,\beta}(\mathbb{R}_+)$, then convolution operator \eqref{eq2.1} is a bounded function $ \forall x\in \R_+$. Moreover, the following inequality holds
	\begin{equation}\label{eq3.12}
	\big\|f \overset{\gamma}{\underset{F_c,F_s,K}{\ast}} g\big\|_{L_{\infty}(\mathbb{R}_+)} \leq 2^{\frac1q}\|f\|_{L_p(\mathbb{R}_+)}\|g\|_{L_q^{0,\beta}(\mathbb{R}_+)},\quad 0<\beta\leq 1.
	\end{equation}

\end{theorem}
\begin{proof}
	Applying  H\"older's inequalities for the pair of conjugate exponents $p$ and $q$, we deduce that 
	
	$$\begin{aligned}
	|(f \overset{\gamma}{\underset{F_c,F_s,K}{\ast}} g)|&\leq \frac14\int_{\mathbb{R}_+^2}|\varphi(x,u,v)||f(u)||g(v)| du dv\\
	&\leq \frac14\left\{\int_{\mathbb{R}_+^2} |\varphi(x,u,v)||f(u)|^p du dv\right\}^{\frac1p}\left\{\int_{\mathbb{R}_+^2}|\varphi(x,u,v)| |g(v)|^q du dv\right\}^{\frac1q}\\
	&=\frac14\left\{\int_{\R_+} |f(u)|^p\left(\int_{\R_+}|\varphi(x,u,v)|dv\right)du\right\}^{\frac1p}
\times\left\{\int_{\R_+} |g(v)|^q\left(\int_{\R_+}|\varphi(x,u,v)|du\right)dv\right\}^{\frac1q}.\end{aligned}$$
Based on \eqref{danhgiaK0}, we have
\begin{equation}\label{eq3.13.1}
\begin{aligned}
|(f \overset{\gamma}{\underset{F_c,F_s,K}{\ast}} g)|&\leq \frac14\left\{\int_{\R_+} |f(u)|^p\left(\int_{\R_+} 4e^{-v} dv\right)du\right\}^{\frac1p}\left\{\int_{\R_+}|g(v)|^q\cdot 8K_0(\beta v) dv\right\}^{\frac1q}\\
&=\frac14\cdot 4^{\left(\frac1p+\frac1q\right)}2^{\frac1q}\|f\|_{L_p(\mathbb{R}_+)}\|g\|_{L_q^{0,\beta}(\mathbb{R}_+)}<\infty.
\end{aligned}
\end{equation}
	Inequality \eqref{eq3.13.1} implies that the convolution operator ($f \overset{\gamma}{\underset{F_c,F_s,K}{\ast}} g$) is a bounded function $\forall x\in \R_+$ and infer the desired conclusion of inequality \eqref{eq3.12}.
\end{proof}
\noindent Overall, through Theorem \ref{thm2.1}, Corollary \ref{cor3.1}, and Theorem \ref{thm3.2} then the characteristic for boundedness of convolution \eqref{eq2.1} on $L_r (\R_+)$  is valid and  well-defined with $r \in [1,\infty]$. Now we will show boundedness on a two-parametric family of Lebesgue spaces $L_s^{\gamma_1,\gamma_2}(\mathbb{R}_+)$ of operator \eqref{eq2.1} as follows.

\begin{theorem}\label{rem3}
	If $p,q$ be real numbers in $(1,\infty)$ such that $\frac1p+\frac1q=1$. For any function $ f\in L_p(\mathbb{R}_+)$ and $g\in L_q^{0,\beta}(\mathbb{R}_+)$, $0<\beta\leq 1$, then convolution \eqref{eq2.1}  is well-defined as
	continuous functions and belonging to $L_s^{\gamma_1,\gamma_2}(\mathbb{R}_+)$ with $s\geq 1$ and two-parametric $\gamma_1 >-1$, $\gamma_2>0$. Moreover
	\begin{equation}\label{eq3.13}
	\big\|f \overset{\gamma}{\underset{F_c,F_s,K}{\ast}} g\big\|_{L_s^{\gamma_1,\gamma_2}(\mathbb{R}_+)} \leq \textup{Const}.\|f\|_{L_p(\mathbb{R}_+)} \|g\|_{L_q^{0,\beta}(\mathbb{R}_+)},
	\end{equation}
	where the upper bound constant in the right-hand side of inequality \eqref{eq3.13} determined by $2^{\frac1q}\gamma_2^{\frac{1-\gamma_1 }{s}}\Gamma^{\frac1s}(\gamma_1 +1)$. Here $L_s^{\gamma_1,\gamma_2}(\mathbb{R}_+)$ is a two-parametric family of Lebesgue spaces defined by $\left\{f(x): \int_{\mathbb{R}_+} |f(x)|^s x^{\gamma_1} e^{-\gamma_2 x} dx<\infty\right\}
	$, 	with the norm
	$
	\|f\|_{L_s^{\gamma_1,\gamma_2}(\mathbb{R}_+)} = \left\{\int_{\mathbb{R}_+} |f(x)|^s x^{\gamma_1} e^{-\gamma_2 x} dx\right\}^{\frac1s}.
	$
\end{theorem}
\begin{proof}
	From inequality \eqref{eq3.13.1}, we infer
	$
	|(f \overset{\gamma}{\underset{F_c,F_s,K}{\ast}} g)(x)|\leq 2^{\frac1q}\|f\|_{L_p(\mathbb{R}_+)}\|g\|_{L_q^{0,\beta}(\mathbb{R}_+)}
	$ is finite. Setting this value as a positive constant $M$.
	According to the formula $3.225$ in \cite{gradshteyn2014Ryzhik}  with $\gamma_1 >-1,\gamma_2 >0$, we deduce that $\int_{\R_+} x^{\gamma_1} e^{-\gamma_2 x} dx=\gamma_2^{(1-\gamma_1)} \Gamma(\gamma_1 +1)$. Therefore
$$
	\int_{\R_+}  x^{\gamma_1} e^{-\gamma_2 x}\big|(f \overset{\gamma}{\underset{F_c,F_s,K}{\ast}}g)(x)\big|^s dx \leq \gamma_2^{(1-\gamma_1)}\Gamma(\gamma_1+1).M^s<\infty,
	$$
	implies that the convolution operator $(f\overset{\gamma}{\underset{F_c,F_s,K}{\ast}} g) \in L_s^{\gamma_1,\gamma_2}(\mathbb{R}_+)$. Moreover, we obtain the following estimation
	$
	\big\|f \overset{\gamma}{\underset{F_c,F_s,K}{\ast}} g\big\|_{L_s^{\gamma_1,\gamma_2}(\mathbb{R}_+)}\leq\underbrace{2^{\frac1q}\bigg(\int_{\R_+} x^{\gamma_1} e^{-\gamma_2 x} dx\bigg)^{\frac1s}}_{\textup{Const}}  \|f\|_{L_p(\mathbb{R}_+)} \|g\|_{L_q^{0,\beta}(\mathbb{R}_+)},
	$ with $s\geq 1$, $\gamma_2 >0$, $\gamma_1>-1$.
	This leads to  $\textup{Const}=2^{\frac1q}\gamma_2^{\frac{1-\gamma_1}{s}}\Gamma^{\frac1s}(\gamma_1 +1)$, where $\Gamma(x)$ is a Euler's gamma-function (refer\cite{Yaku94Luchko}).
\end{proof}
\subsection{Estimation on weighted space $L_p(\mathbb{R}_+,\rho)$}
By considering the $L_p$ norms in more naturally determined weighted spaces. Using the general theory of reproducing kernels, in \cite{Saitoh2000} (also see\cite{Saitoh2016}) Saitoh gave a new inequality for the Fourier convolution in weighted  $L_{p}(\mathbb{R},|\rho_j|)$ Lebesgue spaces as follows
$$\big\lv\big( (F_1 \rho_1) \underset{F}{*}(F_2 \rho_2)\big) \cdot (\rho_1 \underset{F}{*} \rho_2)^{\frac{1}{p}-1}\big\lv_{L_p (\R)} \leq \big\lv F_1 \big\lv_{L_p (\R,|\rho_1|)} \big\lv F_2\big\lv_{L_p (\R,|\rho_2|)},\quad p>1,$$ where $\rho_j $ are non-vanishing functions, $F_j \in L_p (\R, |\rho_j|),\, j=1,2$. Here, the norm of $F_j$ in the weighted space $L_{p}(\mathbb{R}, \rho_j)$ is understood as $
\big\lv F_j  \big\lv_{L_{p}(\mathbb{R}, \rho_j)}=\big\{\int_{\R}|F_j(x)|^{p} \rho_j(x) \mathrm{d} x\big\}^{\frac{1}{p}}.$ This type of inequality is very convenient as many applications require the \textquotedblleft same\textquotedblright $L_p$ norms. It is worth noting that if $f, g$ are functions belonging to $L_2 (\R)$, then Saitoh's inequality for Fourier convolution is still true for, while the same does not happen with Young's inequality. Following this approach, together with using H\"older's inequality, and Fubini's theorem, we establish another result in weighted space $L_{p}(\mathbb{R_+}, \rho_j)$ for convolution operator \eqref{eq2.1}. Some techniques used in the proof of our theorem come from \cite{Tuan2023VKT,2023JmathSci}, and we follow closely the strategy of these results.
\begin{theorem}[Saitoh's type inequality]\label{thm4.1}
	Suppose that $\rho_1, \rho_2$  are non-vanishing  positive functions such that convolution $(\rho_1 \overset{\gamma}{\underset{F_c,F_s,K}{\ast}} \rho_2)$ given by \eqref{eq2.1} is well-defined. For any functions $F_1 \in L_p (\R_+, \rho_1)$ and $F_2 \in L_p (\R_+, \rho_2)$ with $p >1$,  the following $L_p (\R_+)$-weighted inequality holds true 
	\begin{equation}\label{eq4.1}
	I:=\big\|(F_1\rho_1 \overset{\gamma}{\underset{F_c,F_s,K}{\ast}} F_2\rho_2)\cdot (\rho_1 \overset{\gamma}{\underset{F_c,F_s,K}{\ast}} \rho_2)^{\frac1p-1}\big\|_{L_p(\mathbb{R}_+)} \leq (2K_0(v))^{\frac1p} \prod_{i=1}^2\|F_i\|_{L_p(\mathbb{R}_+,\rho_i)}.
	\end{equation}
\end{theorem}
\begin{proof}
	From Definition \eqref{eq2.1}, we obtain
	\begin{equation}\label{eq4.2}
	\begin{aligned}
	I^p&=\big\|(F_1\rho_1 \overset{\gamma}{\underset{F_c,F_s,K}{\ast}} F_2\rho_2)(\rho_1 \overset{\gamma}{\underset{F_c,F_s,K}{\ast}} \rho_2)^{\frac1p-1}\big\|^p_{L_p(\mathbb{R}_+)} \\
	&=\int_{\mathbb{R}_+} \left\{\left|\frac14\int_{\mathbb{R}_+^2} \varphi(x,u,v)(F_1\rho_1)(u)(F_2\rho_2)(v) du dv\right|^p
	\times\left|\frac14\int_{\mathbb{R}_+^2} \varphi(x,u,v)\rho_1(u)\rho_2(v) du dv\right|^{1-p}\right\} dx\\
	&\leq \frac14\int_{\mathbb{R}_+}\left\{\left(\int_{\mathbb{R}_+^2}|\varphi(x,u,v)| |F_1\rho_1(u)| |(F_2\rho_2)(v)| du dv\right)^p
	\times \left(\int_{\mathbb{R}_+^2} |\varphi(x,u,v)| \rho_1(u) \rho_2(v) du dv\right)^{1-p}\right\} dx,
	\end{aligned}\end{equation}
	here, $\varphi(x,u,v)$ is defined by \eqref{eq2.2}.
	Applying the H\"older inequality for the conjugate pair $p$, $q$ we have
	
	\begin{equation}\label{eq4.4}
	\begin{aligned}
	&\int_{\mathbb{R}_+^2} |\varphi(x,u,v)| |F_1(u)| \rho_1(u) |F_2(v)| \rho_2(v) du dv\\
	&\leq \left\{\int_{\mathbb{R}_+^2} |\varphi(x,u,v)| |F_1(u)|^p \rho_1(u) |F_2(v)|^p \rho_2(v) du dv\right\}^{\frac1p}
	\times \left\{\int_{\mathbb{R}_+^2} |\varphi(x,u,v)| \rho_1(u) \rho_2(v) du dv\right\}^{\frac1q}.
	\end{aligned}
	\end{equation}
	Combining \eqref{eq4.4} and \eqref{eq4.2}, we obtain
	\begin{align*}
	I^p&\leq \frac14\int_{\mathbb{R}_+} \left\{\left(\int_{\mathbb{R}_+^2} |\varphi(x,u,v)| |F_1(u)|^p \rho_1(u) |F_2(v)|^p \rho_2(v) du dv\right)\left(\int_{\mathbb{R}_+^2} |\varphi(x,u,v)| \rho_1(u)\rho_2(v) du dv\right)^{\frac{p}{q}+1-p}\right\} dx.
	\end{align*}
	Since $p$ and $q$ are a conjugate pair ($\frac{1}{p}+\frac{1}{q}=1$), it implies that $\frac{p}{q}+1-p=0$, by using Fubini's theorem together with \eqref{danhgiaK0} (evaluation depends on the variable $x$), we infer that 
	\begin{align*}
	I^p&\leq \frac14 \int_{\mathbb{R_+}^3} |\varphi(x,u,v)| |F_1(u)|^p \rho_1(u) |F_2(v)|^p \rho_2(v) du dv dx\\
	&=\frac14\left(\int_{\mathbb{R}_+} |\varphi(x,u,v)| dx\right)\left(\int_{\mathbb{R}_+} |F_1(u)|^p \rho_1(u) du\right)\left(\int_{\mathbb{R}_+} |F_2(v)|^p\rho_2(v) dv\right)\\
	&\leq \frac14 8 K_0(v) \|F_1\|^p_{L_p(\mathbb{R}_+,\rho_1)} \|F_2\|^p_{L_p(\mathbb{R}_+,\rho_2)}.
	\end{align*}
\end{proof}

In case one of functions $\rho_1(x)$, $\rho_2(x)$ is homogenous 1, for instance $\rho_1(x)\equiv 1$ for all $ x\in\mathbb{R}_+$, and $0<\rho_2\in L_1(\mathbb{R}_+)$, then we have
\begin{align*}
|(1 \overset{\gamma}{\underset{F_c,F_s,K}{\ast}} \rho_2)(x)| &\leq \frac14 \int_{\mathbb{R_+}^2} |\varphi(x,u,v)| \rho_2(v) du dv
=\frac14\left(\int_{\mathbb{R}_+} |\varphi(x,u,v)| du\right)\left(\int_{\mathbb{R}_+} \rho_2(v) dv\right)\\ &\leq 2K_0(v)\|\rho_2\|_{L_1(\mathbb{R}_+)}<\infty.
\end{align*}
This means that $(1 \overset{\gamma}{\underset{F_c,F_s,K}{\ast}} \rho_2)(x)$ is well-defined and therefore 
$|(1 \overset{\gamma}{\underset{F_c,F_s,K}{\ast}} \rho_2)(x)|^{1-\frac1p} \leq \{2K_0(v)\}^{1-\frac1p}\|\rho_2\|^{1-\frac1p}_{L_1(\mathbb{R}_+)}.$ Combining with Theorem \ref{thm4.1}, we arrive at the following corollary.

\begin{corollary}\label{cor4.1}
	Let $\rho_2 $ be a  positive function belonging to $L_1 (\R_+)$. If  $F_1,F_2$ are functions belonging to  $L_p (\R_+)$ and $L_p (\R_+,\rho_2)$, respectively, with $p>1$, then  the following estimate holds true	
	\begin{equation}\label{eq4.6}
	\| F_1 \overset{\gamma}{\underset{F_c,F_s,K}{\ast}} F_2 \rho_2\|_{L_p(\mathbb{R}_+)} \leq 2K_0(v) \|\rho_2\|^{1-\frac1p}_{L_1(\mathbb{R}_+)} \|F_1\|_{L_p(\mathbb{R}_+)} \|F_2\|_{L_p(\mathbb{R}_+,\rho_2)}.
	\end{equation}
\end{corollary}
\noindent For example, choose $\rho_1(x)\equiv 1$, $\forall x\in \mathbb{R}_+$ and $\rho_2(x)=e^{-x}\in L_1(\mathbb{R}_+)$:
$
|(1 \overset{\gamma}{\underset{F_c,F_s,K}{\ast}} e^{-x})(t)| \leq 2K_0(v)
$ with $v>0$,
and
$
\|F_1 \overset{\gamma}{\underset{F_c,F_s,K}{\ast}} F_2 \rho_2\|_{L_p(\mathbb{R}_+)} \leq 2K_0(v)\|F_1\|_{L_p(\mathbb{R}_+)} \|F_2\|_{L_p(\mathbb{R}_+, e^{-x})}.
$
\begin{remark}\label{rem4}
\textup{	The Young type estimation \eqref{eq3.10} is not valid in the typical space $L_2(\mathbb{R}_+)$, while the estimations \eqref{eq4.1} and \eqref{eq4.6} are hold on $L_2(\mathbb{R}_+)$ space.}
\end{remark}

\section{$L_1$-solution for classes of convolution integro-differential equations}
This section will be devoted to a class of the first and second kind of convolution integral equations
related to \eqref{eq2.1}. We will establish conditions that will
guarantee the existence and uniqueness of solutions in a closed form for these equations. Similar
questions for different convolutions were considered in \cite{Yaku94Luchko,Yakubovich1996index}. 
 The solutions are presented in closed form via the convolution \eqref{eq2.1}. We obtain the boundedness of solutions on the space $L_p(\mathbb{R}_+)$, $p\geq 1$.

\subsection{Second kind of convolution integro-differential equation }
We consider the following equation  $f(x) + \left(1 - d^2\textfractionsolidus dx^2\right) \frac{1}{\sqrt{2\pi}}\int_{\R_+} f(u)[g(|x-u|) + g(x+u)] du = h(x),\ x>0.$
Taking into account the symmetric properties of the convolution kernel of  $(F_c)$ transform \eqref{eq1.6}, then the above equation can be written as follow
\begin{equation}\label{eq5.1}
f(x) + \left(1 - \frac{d^2}{dx^2}\right) (f \underset{F_c}{\ast} g) (x) = h(x),\quad x>0.
\end{equation}
Here $f(x)$ is an unknown function need to find. We will give the conditions for solvability in $L_1(\R_+)$ of equation \eqref{eq5.1} for the case $g$, $h$ are given functions, defined by $h=(\varphi \overset{\gamma}{\underset{F_c,F_s,K}{\ast}} \xi)$ and $g(x)=(\mathrm{sech}\, t \underset{F_c}{\ast} g_1(t))(x)$, where $g_1,\varphi \in L_1(\mathbb{R}_+)$, $\xi \in L_1^{0,\beta}(\mathbb{R}_+)$ with $\beta \in (0,1)$.

\begin{theorem}\label{thm5.1}
	Let $g_1$, $\varphi$ and $\xi$ are given functions such that $g_1, \varphi\in L_1(\mathbb{R}_+)$ and $\xi\in L_1^{0,\beta}(\mathbb{R}_+)$. Then
	for the solvability of \eqref{eq5.1} in $L_1 (\R_+)$, those are the sufficient conditions include 
	$
	F_c\big(\mathrm{sech}^3t \underset{F_c}{\ast} g_1(t)\big) (y) \ne 0$ and that $(1+y^2)F_c\big(\mathrm{sech}^3t \underset{F_c}{\ast} g_1(t)\big)(y)$ is finite,  for any $y>0$. Moreover, its solution $f(x)$ is unique and represented by formula 
	$
	f(x) = \big(\varphi \overset{\gamma}{\underset{F_c,F_s,K}{\ast}} \xi \big)(x) -\big(\ell \underset{F_c}{\ast}\big(\varphi \overset{\gamma}{\underset{F_c,F_s,K}{\ast}} \xi\big) \big)(x)
	$ almost everywhere in $\R_+$, where $\ell \in L_1(\mathbb{R}_+)$ is defined via
	$
	(F_c\ell)(y)=\frac{2F_c\big(\mathrm{sech}^3 t\underset{F_c}{\ast} g_1\big)(y)}{1+F_c\left(\mathrm{sech}^3 t\underset{F_c}{\ast} g_1\right)(y)}.
	$
 Finally, the following $L_1$-norm estimation hold
	$
	\|f\|_{L_1(\mathbb{R}_+)} \leq 2\|\varphi\|_{L_1(\mathbb{R}_+)} \|\xi\|_{L_1^{0,\beta}(\mathbb{R}_+)}\big(1 + 2\sqrt{\frac{2}{\pi}} \|\ell\|_{L_1(\mathbb{R}_+)}\big)
	$.
\end{theorem}

To prove Theorem \ref{thm5.1}, first, we need the following auxiliary lemma.
\begin{lemma}\label{lem5.1} Let $f \in L_1 (\R_+)$ and
	$g$ is an $L_1$-Lebesgue integrable function on $\R_+$ satisfying condition $(1+y^2)(F_cg)(y)$ is finite. For any $y>0$, we have the following assertion
	\begin{equation}\label{eq5.2}
	F_c\left[\left(1-\frac{d^2}{dx^2}\right)\big(f\underset{F_c}{\ast} g\big) (x) \right](y) = (1+y^2)(F_cf)(y)(F_cg)(y).
	\end{equation}
\end{lemma}
\begin{proof}
	Using the Parseval equality for the convolution \eqref{eq1.6} we obtain
	\begin{equation}\label{eq5.3}
	\big(f \underset{F_c}{\ast} g\big)(x) = \sqrt{\frac{2}{\pi}} \int_{\mathbb{R}_+} (F_cf)(y)(F_cg)(y) \cos (xy) dy,\quad x>0.
	\end{equation}
	Under the $f,g\in L_1(\mathbb{R}_+)$ assumption, we obviously have both
 $(F_cf)(y)$ and $(F_cg)(y)$ belong to $C_0(\mathbb{R}_+)$ \cite{Titchmarsh1986}. Therefore 
	$
	\|(F_cf)(y)| \leq \|f\|_{L_1(\mathbb{R}_+)}$ is finite and  $\|(F_cg)(y)| \leq \|g\|_{L_1(\mathbb{R}_+)} 
	$ is finite,
	these imply that the product of $(F_cf)(y)(F_cg)(y)$ is a  bounded function on $\R_+$. Then, the integral in right-hand side of \eqref{eq5.3} is absolute convergent, so  we can change the differential and the integral order as follows
	\begin{align*}
	\frac{d^2}{dx^2}( f\underset{F_c}{\ast} g)(x)&= \sqrt{\frac{2}{\pi}} \frac{d^2}{dx^2} \int_{\mathbb{R}_+} (F_cf)(y) (F_cg)(y) \cos (xy) dy
	=\sqrt{\frac{2}{\pi}}\int_{\mathbb{R}_+}(F_cf)(y)(F_cg)(y)\left(\frac{d^2}{dx^2}\cos (xy)\right) dy\\
	&=\sqrt{\frac{2}{\pi}}\int_{\mathbb{R}_+}(-y^2)(F_cf)(y)(F_cg)(y)\cos(xy)dy
=F_c(-y^2(F_cf)(y)(F_cg)(y))(x).
	\end{align*}
	Therefore,
	\begin{equation}\label{eq5.4}
	\left(1-\frac{d^2}{dx^2}\right)\big(f\underset{F_c}{\ast} g\big)(x)=F_c((1+y^2) (F_cf)(y) (F_cg)(y))(x).
	\end{equation}
	Furthermore, $(1+y^2)(F_cg)(y)< \infty$, together with $|(F_cf)(y)| \leq \|f\|_{L_1(\mathbb{R}_+)}$. We deduce $(1+y^2)(F_cf)(F_cg)(y)$ belongs to $L_1(\mathbb{R}_+)$. Applying the Fourier cosine transform on both sides of \eqref{eq5.4},we come to the conclusion of lemma.
\end{proof}
It should be emphasized that the condition $(1+y^2)(F_cg)(y)< \infty$ in Lemma \ref{lem5.1} is always existence with the given assumptions. A straightforward instance shows this if we choose $g=\sqrt{\frac{\pi}2} e^{-x}$. It is easy to check this function belongs to $L_1(\mathbb{R}_+)$, then $(F_cg)(y)=\frac{1}{1+y^2}$. Therefore $(1+y^2)F_c\left(\sqrt{\frac{\pi}2}e^{-x}\right)(y)=1$. 
\begin{proof}[Proof of Theorem \ref{thm5.1}]
	Applying $(F_c)$-transform on both sides of equation \eqref{eq5.1}, for any $y>0$, we have $(F_cf)(y) + (1+y^2)(F_cf)(y)(F_cg)(y) = (F_ch)(y)$. Using \eqref{eq5.2} we obtain \begin{equation}\label{eq5.5}
	(F_cf)(y)\left[1+ (1+y^2) F_c\big(\mathrm{sech}\, t \underset{F_c}{\ast} g_1(t)\big)(y)\right] = (F_ch)(y).
	\end{equation}
	According to \cite{bateman1954}, using formulas 1.9.1, and 1.9.4 page 30, we have respectively	
	$
	(F_c\mathrm{sech}\, t)(y) = \sqrt{\frac{\pi}{2}} \mathrm{sech}\,\frac{\pi y}{2},$ and $ \frac{\sqrt{2\pi}}{4}(1+y^2)\mathrm{sech}\,\frac{\pi y}{2}=F_c(\mathrm{sech}^3 t)(y).
	$
This leads to
	\begin{align*}
	(1+y^2)F_c\big(\mathrm{sech}\, t \underset{F_c}{\ast} g_1(t)\big)(y) &= (1+y^2)F_c(\mathrm{sech}\,t) (y)(F_cg_1)(y)
	= (1+y^2)\sqrt{\frac{\pi}{2}}\mathrm{sech}\,\frac{\pi y}{2} (F_cg_1)(y)\\
	&=2F_c(\mathrm{sech}^3 t)(y)(F_cg_1)(y)
	= 2F_c\big(\mathrm{sech}^3t\underset{F_c}{\ast} g_1(t)\big)(y).
	\end{align*}
	From \eqref{eq5.5} and under the given condition, we have
	$
	(F_cf)(y)=\frac{1}{1+2F_c\big(\mathrm{sech}^3t \underset{F_c}{\ast} g_1(t)\big)(y)}\cdot (F_ch)(y)
	$
	or equivalent to
	$
	(F_cf)(y)=\bigg(1-\frac{2F_c\left(\mathrm{sech}^3 t\underset{F_c}{\ast} g_1(t)\right)(y)}{1+2\left(F_c\mathrm{sech}^3t \underset{F_c}{\ast} g_1(t)\right)(y)}\bigg)(F_ch)(y).
	$
	Apply the Wiener--L\'evy's theorem \cite{NaimarkMA1972} for the Fourier cosine $(F_c)$ transform, there exists a function $\ell \in L_1(\mathbb{R}_+)$ such that
	$
	(F_c\ell)(y) = \frac{2F_c\left(\mathrm{sech}^3t \underset{F_c}{\ast} g_1\right)(y)}{1+ 2\left(F_c\mathrm{sech}^3 t\underset{F_c}{\ast} g_1\right)(y)}, 
	$ for any $y>0$.
	This yields $(F_cf)(y)=(F_ch)(y)-(F_c\ell)(y)(F_ch)(y)$, or equivalent to
	$$(F_cf)(y)=F_c\big(\varphi \overset{\gamma}{\underset{F_c,F_s,K}{\ast}} \xi\big)(y) - F_c\left(\ell \underset{F_c}{\ast} \big(\varphi \overset{\gamma}{\underset{F_c,F_s,K}{\ast}} \xi\big)\right)(y).$$
	Therefore
	$
	f(x)=\big(\varphi \overset{\gamma}{\underset{F_c,F_s,K}{\ast}} \xi\big)(x) - \left(\ell \underset{F_c}{\ast} \big( \varphi \overset{\gamma}{\underset{F_c,F_s,K}{\ast}} \xi \big)\right)(x)$ almost everywhere on $\mathbb{R}_+.
	$
	Now we need to point out $f(x) \in L_1 (\R_+)$. Since $\varphi \in L_1(\mathbb{R}_+)$ and  $\xi \in L_1^{0,\beta}(\mathbb{R}_+)$, by Theorem \ref{thm2.1}, we infer that $(\varphi \overset{\gamma}{\underset{F_c,F_s,K}{\ast}} \xi) $ well-defined as a continuous function and belong to $L_1(\mathbb{R}_+)$. Moreover, we have $\ell \in L_1(\mathbb{R}_+)$ due to Wiener--L\'evy's theorem, which implies that $\ell \underset{F_c}{\ast}(\varphi \overset{\gamma}{\underset{F_c,F_s,K}{\ast}} \xi) \in L_1(\mathbb{R}_+)$. The combination of these is sufficient to show there exists the only existence of $f(x)$ belonging to $L_1(\mathbb{R}_+)$ because convolutions abbreviated as ${\underset{F_c,F_s, K}{\ast}}$ and ${\underset{F_c}{\ast}}$ uniquely determined. By applying inequalities \eqref{eq1.8} and \eqref{eq2.3} we obtain the estimtaion on $L_1(\mathbb{R}_+)$ as follows
	\begin{align*}
	\|f\|_{L_1(\mathbb{R}_+)} &\leq \big\|\varphi \overset{\gamma}{\underset{F_c,F_s,K}{\ast}} \xi\big\|_{L_1(\mathbb{R}_+)}\left(1+ 2\sqrt{\frac{2}{\pi}}\|\ell\|_{L_1(\mathbb{R}_+)}\right)
	\leq 2\|\varphi\|_{L_1(\mathbb{R}_+)} \|\xi\|_{L_1^{0,\beta}(\mathbb{R}_+)}\left(1+2\sqrt{\frac{2}{\pi}}\|\ell\|_{L_1(\mathbb{R}_+)}\right)
	\end{align*}
\end{proof}

\subsection{Integro-differential equation of first type} This remainder will deal with equations of the form\\ $ \left(1 - d^2\textfractionsolidus dx^2\right) \frac{1}{\sqrt{2\pi}}\int_{\R_+} f(u)[g(|x-u|) + g(x+u)] du = h(x).$ Taking into account the symmetric properties of the convolution kernel of  $(F_c)$ transform \eqref{eq1.6}, which can then be written as convolutional equation
\begin{equation}\label{eq5.10}
\left(1 - d^2\textfractionsolidus dx^2\right)\big(f \underset{F_c}{\ast} g\big)(x) = h(x),
\end{equation}
where $f(x)$ is an unknown function need to find. We show the conditions for solvability in $L_1(\R_+)$ of convolution equation \eqref{eq5.10} for the case $g$ being a given functions in $L_1 (\R_+)$, $h:=(\psi \overset{\gamma}{\underset{F_c,F_s,K}{\ast}} \xi)(x),$ where $\xi \in L_1^{0,\beta}(\mathbb{R}_+)$, $\beta \in (0,1)$. Function $\psi$ is determined via by $\psi:=(\varphi \underset{F_s,F_c}{\ast} g)(x)$ with $\varphi$ being a given functions in $L_1 (\R_+)$. Here $\underset{F_s,F_c}{\ast}$ is denoted by the generalized convolution of two functions for the Fourier sine and Fourier cosine transforms \cite{Sneddon1972} defined by
	$(\varphi \underset{F_s,F_c}{\ast} g)(x) := \frac{1}{\sqrt{2\pi}} \int_{\R_+} \varphi(u)[g(|x-u|) - g(x+u)] du,\ x>0.$
	With assumption $\varphi,g \in L_1(\mathbb{R}_+)$ then $(\varphi \underset{F_s,F_c}{\ast} g)$ belongs to $L_1(\mathbb{R}_+)$ and the following factorization equality holds
$	F_s(\varphi \underset{F_s,F_c}{\ast} g)(y) = (F_s\varphi)(y)(F_cg)(y)$ for all $y>0$ (refer \cite{Sneddon1972}). We easily get the following proposition.
\begin{proposition}\label{menhde4.1}
	Let $\varphi, g\in L_1(\mathbb{R}_+)$, then we have estimation 
	$\|\varphi \underset{F_s,F_c}{\ast} g\|_{L_1(\mathbb{R}_+)} \leq 2\sqrt{\frac{2}{\pi}} \|\varphi\|_{L_1(\mathbb{R}_+)} \|g\|_{L_1(\mathbb{R}_+)}.$	
For $p, q, r\in (1,\infty)$ such that $\frac1p+\frac1q=1+\frac1r$. For any functions $\forall \varphi \in L_p(\mathbb{R}_+)$, $g\in L_q(\mathbb{R}_+)$ then $(\varphi \underset{F_s,F_c}{\ast} g)$ well-defined as a continuous function belong to $L_r(\mathbb{R}_+)$. Moreover, the following estimation holds
$\|\varphi \underset{F_s,F_c}{\ast} g\|_{L_r(\mathbb{R}_+)} \leq 2\sqrt{\frac{2}{\pi}}\|\ell\|_{L_p(\mathbb{R}_+)} \|h\|_{L_q(\mathbb{R}_+)}.$
\end{proposition}

\begin{theorem}\label{thm5.2} With the assumptions in the description of equation \eqref{eq5.10}. 
	Suppose that $g$ is an $L_1$-Lebesgue integrable function over $\R_+$ such that $(F_cg)(y)\ne 0$ and $(1+y^2)(F_cg)(y)$ is finite, for any $y>0$. Then equation \eqref{eq5.10} has the unique solution in $L_1(\mathbb{R}_+)$ which can be present by
	$
	f(x) = \big((\varphi \underset{F_s,F_c}{\ast} \sqrt{\frac{2}{\pi}} e^{-t}) \overset{\gamma}{\underset{F_c,F_s,K}{\ast}} \xi\big)(x)
	$ almost everywhere on $\R_+$.
	Moreover, the following estimation holds
	$
	\|f\|_{L_1(\mathbb{R}_+)} \leq 4\|\varphi\|_{L_1(\mathbb{R}_+)} \|\xi\|_{L_1(\mathbb{R}_+)}.
	$
\end{theorem}
\begin{proof}
	Applying the Fourier cosine $(F_c)$ on both sides of equation \eqref{eq5.10}, and combining with the results \eqref{eq5.2}, and \eqref{eq2.5}  we obtain
	\begin{equation}\label{eq5.11}
	(1+y^2)(F_cf)(y)(F_cg)(y) = \sin y(F_s\psi)(y)K[\xi](y),\quad y>0.
	\end{equation}
	Due to the way the $\psi(x)$ function is defined by $\psi(x)=(\varphi \underset{F_s,F_c}{\ast} g)(x)$, based on the factorization equality of $\underset{F_s,F_c}{\ast}$ generalized convolution \cite{Sneddon1972}, we deduce that
	$
	(F_s\psi)(y)=F_s(\varphi \underset{F_s,F_c}{\ast} g)(y)=(F_s\varphi)(y)(F_cg)(y)
	$ holds for any $y>0$.
	Coupling \eqref{eq5.11} and condition $(F_cg)y\ne 0$, we obtain
	$(1+y^2)(F_cf)(y) = \sin y(F_s\varphi)(y) K[\xi](y)$, this equivalent to
	$(F_cf)(y)=\frac{1}{1+y^2} \sin y(F_s\varphi)(y)K[\xi](y).$
Moreover
	$F_c\left(\sqrt{\frac{\pi}{2}}e^{-t}\right)(y)=\frac{1}{1+y^2}$. Therefore
	$$
	(F_cf)(y)=\sin y F_c\left(\sqrt{\frac{\pi}{2}}e^{-t}\right)(y) (F_s\varphi)(y) K[\xi](y).
	$$
	Again using factorization equality of $(F_s,F_c)$-generalized convolution together with \eqref{eq2.5}, we infer 
	$$
	(F_cf)(y)=\sin y F_s\left(\varphi \underset{F_s,F_c}{\ast} \sqrt{\frac{\pi}{2}} e^{-t}\right)(y) K[\xi](y)
	=F_c\left[\left(\varphi \underset{F_s,F_c}{\ast} \sqrt{\frac{\pi}{2}} e^{-t}\right)\overset{\gamma}{\underset{F_c,F_s,K}{\ast}} \xi\right](y).$$ This means that
	$f(x)=\left(\big(\varphi \underset{F_s,F_c}{\ast} \sqrt{\frac{\pi}{2}} e^{-t}\big) \overset{\gamma}{\underset{F_c,F_s,K}{\ast}} \xi\right)(x)$ almost everywhere on $\R_+$.
	Since $\varphi \in L_1(\mathbb{R}_+)$, then obviously that $(\varphi \underset{F_s,F_c}{\ast}\sqrt{\frac{\pi}{2}} e^{-t})$ belongs to $L_1(\mathbb{R}_+)$ (refer \cite{Sneddon1972}). Furthermore $\xi \in L_1^{0,\beta}(\mathbb{R}_+)$ we directly deduce $f(x)$ belongs to $L_1(\mathbb{R}_+)$ based on Theorem \ref{thm2.1}. For any $\varphi, \xi$ are given function in $L_1(\mathbb{R}_+)$, using  \eqref{eq2.3} and estimations in Proposition \ref{menhde4.1}, we obtain
	\begin{align*}
	\|f\|_{L_1(\mathbb{R}_+)} &\leq 2\bigg\|\varphi \underset{F_s,F_c}{\ast} \sqrt{\frac{\pi}{2}} e^{-t}\bigg\|_{L_1(\mathbb{R}_+)} \|\xi\|_{L_1(\mathbb{R}_+)}
\leq 4\sqrt{\frac{2}{\pi}} \|\varphi\|_{L_1(\mathbb{R}_+)} \left\|\sqrt{\frac{\pi}{2}} e^{-t}\right\|_{L_1(\mathbb{R}_+)} \|\xi\|_{L_1(\mathbb{R}_+)}\\
	 &\leq 4\|\varphi\|_{L_1(\mathbb{R}_+)} \|\xi\|_{L_1(\mathbb{R}_+)},
	\end{align*}	
\end{proof}
We give another estimate for the solution through Proposition \ref{menhde4.1} and inequality \eqref{eq3.10} as follow
\begin{remark}\label{rem5}
\textup{Let $p,q,r$ are real numbers in open interval $(1,\infty)$ such that $\frac1p+\frac1r=1+\frac1q$, then we have\\
	$
	\|\varphi \underset{F_s,F_c}{\ast} \sqrt{\frac{\pi}{2}} e^{-t}\|_{L_q}\leq 2\left(\frac1r\right)^{\frac1r} \|\varphi\|_{L_p(\mathbb{R}_+)}.
	$
If $\frac1p+\frac1q+\frac1r=2+\frac1s$, $s>1$. Then we obtain the solution's estimation for the problem \eqref{eq5.10} as follow
	$
	\|f\|_{L_s(\mathbb{R}_+)} \leq 2^{1+\frac1q} \left(\frac1r\right)^{\frac1r}\|\varphi\|_{L_p(\mathbb{R}_+)}\|\xi\|_{L_q^{0,\beta}(\mathbb{R}_+)},
	$
	$\forall f\in L_s(\mathbb{R}_+)$, $\varphi\in L_p(\mathbb{R}_+)$ and  $\xi \in L_q^{0,\beta}(\mathbb{R}_+)$ with $\beta \in (0,1)$.}
\end{remark}

\vskip 0.3cm
\noindent \textbf{Disclosure statement}\\
\noindent No potential conflict of interest was reported by the authors.
\vskip 0.3cm
\noindent \textbf{Funding}\\
\noindent This research received no specific grant from any funding agency.
\vskip 0.3cm
\noindent \textbf{ORCID}\\
\noindent \textit{Trinh Tuan} {\color{blue} \url{https://orcid.org/0000-0002-0376-0238}}

\end{document}